\documentclass[11pt]{article}
\usepackage{amsmath}
\usepackage{amsfonts}
\usepackage{amssymb}
\usepackage{amsthm}
\usepackage{qtree}
\usepackage{xcolor}
\usepackage{dsfont}

 \allowdisplaybreaks[4]

\newcommand{\N}{{\mathbb N}}
\newcommand{\M}{{\mathcal M}}

\newcommand{\R}{{\mathbb R}}
\newcommand{\LL}{{\mathcal L}}

\newtheorem{theorem}{Theorem}
\newtheorem{corollary}{Corollary}

\DeclareMathOperator{\ri}{ri}
\DeclareMathOperator{\var}{var}
\DeclareMathOperator{\diam}{diam}

\newtheorem{lemma}{Lemma}
\newtheorem{remark}{Remark}
\author{GUAN-ZHONG MA\thanks{mgz09@mails.tsinghua.edu.cn} \quad YAO XIAO \thanks{yaoxiao0710107@163.com}
\\Department of Mathematical Sciences, \\Tsinghua University,
Beijing, China}
\date{}
\title{Higher dimensional multifractal analysis of non-uniformly hyperbolic systems}

\begin{document}
\maketitle
\begin{abstract}

Johansson, Jordan, \"Oberg and Pollicott ( Israel J. Math.(2010)) has studied the multifractal analysis of a class of one-dimensional non-uniformly hyperbolic systems,
by introducing some new techniques, we extend the results to the case of high dimension.

\smallskip



\end{abstract}

\section{Introduction}
In this note we  extend  recent work on multifractal analysis of non-uniformly hyperbolic system, conducted by Johansson, Jordan, \"Oberg and Pollicott \cite{JJOP2010} to higher dimension. Although the results are quite similar, our methods are different from theirs in several aspects. We will explain it more precisely after we present the main results.

We start with an introduction about the basic settings.
Let $T:\bigcup_{i=1}^m I_i\rightarrow [0,1]$ be a piecewise $C^{1}$ map satisfies the following condition:
\begin{itemize}
\item $I_i\subset [0,1], i=1,\cdots,m$  such that $I_{i}$ and $I_{j}$ does not overlap  for $i\ne j.$
\item $T|_{ I_j}:I_j\rightarrow [0,1]$ is onto and $C^{1}$ map, for all $1\leq j\leq m.$ There is a unique $x_j\in I_j$ such that $T(x_j)=x_j.$
\item   $T'(x)>1$ for $x\not\in\{x_1,\cdots, x_m\}$.
\end{itemize}
We remark that since the map $T$ is $C^1$, we have $T^\prime(x_j)\ge 1$ for $j=1,\cdots,m.$ If for some $j$,
   $T'(x_j)=1$, we call   $x_j$ a {\it parabolic} fixed point.

Define the attractor of $T$ as
$$
  \Lambda=\{x\in \bigcup_{j=1}^m I_j | T^{n}(x)\in [0,1], \forall n\geq0\}.
$$
It is well known  that $\Lambda$ is  invariant under $T$ and we get a  dynamic system
$T: \Lambda \rightarrow \Lambda.$

The above system has a symbolic coding which can be defined as follows.
Let  $T_{i}$ be the inverse map of $T|_{I_{i}}: I_{i}\rightarrow [0,1]$ for $i=1,\cdots,m$.
  Let $\mathcal{A}=\{1\dots, m\}$ and $\Sigma=\mathcal{A}^{\mathbb{N}}$. There is a shift map $\sigma:\Sigma\to\Sigma$ defined by $\sigma((\omega_n)_{n\ge 1})=(\omega_n)_{n\ge 2}$. Define a projection $\Pi: \Sigma \to [0,1]$ as
\[
\Pi(\omega)=\lim_{n\rightarrow\infty}T_{\omega_1}\circ T_{\omega_2}\circ \dots \circ T_{\omega_n}([0,1]).
\]
Then $\Pi(\Sigma)=\Lambda$ and moreover
%\[\begin{array}{cccc}
%\Sigma &
%\stackrel{\sigma}{\longrightarrow} &
%\Sigma  \\
%\Big\downarrow \vcenter{%
%\rlap{$\scriptstyle{\Pi}$}}& &
%\Big\downarrow\vcenter{%
%\rlap{$\scriptstyle{\Pi}$}}\\
%\Lambda &\stackrel{T} {\longrightarrow}&
%\Lambda.
%\end{array}\]
$$
 \Pi\circ\sigma(\omega)=T\circ\Pi(\omega).
$$

%Obviously, $\Pi$ is a bijection except for at most countable points.

Given $F \in C(\Lambda,\R^{d})$ and $\alpha\in\R^d$ one can define the {\it level set} as
\[
\Lambda_{\alpha}=\left\{x\in\Lambda:\lim_{n\rightarrow\infty}\frac{1}{n}\sum_{j=0}^{n-1}F(T^{j}x)=\alpha\right\}.
\]
Not every $\Lambda_\alpha$ is nonempty. Indeed it is also classical that $\Lambda_\alpha\ne \emptyset$ if and only if
$$
\alpha\in \mathcal{L}_{F}:=\left\{\int F d\mu : \mu \in \M(\Lambda,T)\right\}
$$
where  $\mathcal{M}(\Lambda,T)$ is the set of all invariant probability measures on $(\Lambda,T).$

The central problem in multifractal analysis is to determine the size of  $\Lambda_\alpha$, especially the Hausdorff dimension of  it.

Similarly given  $f \in C(\Sigma,\R^{d})$ and $\alpha\in\R^d$ one can define the {\it level set} as
\[X_{\alpha}=\left\{\omega\in\Sigma:\lim_{n\rightarrow\infty}\frac{1}{n}\sum_{j=0}^{n-1}f(\sigma^{j}\omega)=\alpha \right\}.\]
If we denote by $\mathcal{M}(\Sigma,\sigma)$  the set of all invariant probability measures on $(\Sigma,\sigma)$ and define
 $$
 \LL_{f}=\left\{\int fd\mu: \mu \in \mathcal{M}(\Sigma,\sigma)\right\}.
 $$
 Then $X_\alpha\neq \emptyset$ if and only if $\alpha\in \LL_f.$

 Two kinds of level set are related in the following way.  Given $F:\Lambda\to\R^d$ continuous. Define $f:=F\circ\Pi$, then $f:\Sigma\to\R^d$ is continuous and
 $$
 \Pi(X_\alpha)=\Lambda_\alpha.
 $$

Define  $g(\omega):=-\log T'_{\omega_1}\Pi(\sigma \omega)$ and let
$$
\tilde{\Sigma}=\{\omega\in\Sigma:\liminf\limits_{n\rightarrow \infty }\frac{1}{n}\sum_{j=0}^{n-1}g(\sigma^j\omega)>0\}.
$$
Let $h(\mu,\sigma)$, $\lambda(\mu,\sigma)$  be  the metrical entropy and Lyapunov exponent of $\mu$. We have the following  theorem:
\begin{theorem}\label{main-1}
Assume $f:\Sigma\to\R^d$ is continuous. Then for $\alpha\in \LL_f$,
$$\dim_{\text{H}}\Pi(X_{\alpha}\cap \tilde{\Sigma})= \underset{\mu\in \mathcal{M}(\Sigma,\sigma)}{\sup}\left\{\;\frac{h(\mu,\sigma)}{\lambda(\mu,\sigma)}\;\left|\;\int f d \mu=\alpha, \lambda(\mu,\sigma)>0\right.\right\}.
$$
\end{theorem}

Now we will use the above theorem to the non-uniformly hyperbolic system.  Consider the system $T: \Lambda\to\Lambda.$
Let $I\subset \{x_1,\cdots, x_m\}$ be the set of parabolic fixed points. Given $F:\Lambda\to\R^d $ continuous and define
$
A={\rm Co}\{F(x): x\in I\},$ which is  the {\it convex hull} of  $\{F(x):x\in I\}$.

\begin{theorem}\label{main-2}
Assume that $(\Lambda,T)$ is a system defined as above. Given  $F$ continuous  and define $A$  as above. If  for any $\epsilon>0$, there exists $\nu\in\M(\Lambda,T)$ with $\lambda(\nu,T)>0$ and
$\frac{h(\nu,T)}{\lambda(\nu,T)}>\dim_{\text{H}}\Lambda-\epsilon$. Then for $\alpha\in \LL_F \setminus A$, we have
$$
\dim_{\text{H}}\Lambda_{\alpha}= \underset{\mu\in \mathcal{M}(\Lambda,T)}{\sup}\left\{\frac{h(\mu,T)}{\lambda(\mu,T)}\left|\int F d \mu=\alpha\right. \right\},
$$
and   for all $\alpha\in A$ we have $\dim_{\text{H}}\Lambda_{\alpha}=\dim_{\text{H}}\Lambda$.
\end{theorem}

 To present the next result, we need several notations from convex analysis. Given $C\subset\R^d$, the {\it affine hull} of $C$ is the smallest affine subspace of $\R^d$  containing $C$ and is denoted by $\text{aff}( C )$.
For a convex set $C$, we define $\ri ( C )$, the {\it relative interior} of $C$ as
$$\ri  ( C ):=\{x\in{\rm aff} ( C ): \exists \epsilon>0, (x+\epsilon B)\cap\text{aff} ( C )\subset C\},$$
where $B=B(0,1)\subset\R^d$ is the unit open ball.

Write $D(\alpha)=\dim_H \Lambda_\alpha$. We have
\begin{corollary}\label{main-3}
Under the assumption of Theorem 2, we have
$D(\alpha)=\dim_{H}\Lambda$ for  $\alpha\in A$.
   $D(\alpha)$ is continuous in ${\rm ri}(\LL_F)\setminus A.$
\end{corollary}

\begin{remark}
{\rm Our result is a generalization of that in \cite{JJOP2010}. Indeed in \cite{JJOP2010}, they deal with the scalar potential $F:\Lambda\to\R$, while here we deal with vector potential $F:\Lambda\to \R^d.$
}
\end{remark}

In the following we will explain that in the higher dimension case   some extra difficulties occur and  the argument given in \cite{JJOP2010} will no longer work. We will also give the idea that how we overcome these difficulties.

\begin{remark}
{\rm
The first difficulty comes from the proof of the lower bound in Theorem \ref{main-1}. In \cite{JJOP2010}, for each $\alpha$ in the interior of $\LL_f=[\alpha_{\min},\alpha_{\max}]$ and each $\mu\in \M(\Sigma,\sigma)$ such that $\int f d\mu=\alpha$ and $\lambda(\mu,\sigma)>0$, they can construct a sequence $n$-level Bernoulli measures $\mu_n$, which is  $\sigma^n$-ergodic, with
$\int A_n f d\mu_n=\alpha$ and $\lambda(\mu_n,\sigma^n)>0$ such that
$$
\frac{h(\mu_n,\sigma^n)}{\lambda(\mu_n,\sigma^n)}\to \frac{h(\mu,\sigma)}{\lambda(\mu,\sigma)}
$$
(see \cite{JJOP2010} Lemma 3).
The ergodicity of $\mu_n$ implies that $\mu_n(X_\alpha\cap \tilde \Sigma)=1$ and $\dim_H\mu_n=\frac{h(\mu_n,\sigma^n)}{\lambda(\mu_n,\sigma^n)}$. Now  it is known that (see \cite{JJOP2010} Lemma 2) the average of $\mu_n$ will give a $\sigma$-invariant and ergodic measure $\nu_n$ and satisfies $\int f d\nu_n=\alpha$ and
$$
\frac{h(\mu_n,\sigma^n)}{\lambda(\mu_n,\sigma^n)}= \frac{h(\nu_n,\sigma)}{\lambda(\nu_n,\sigma)}.
$$
Then the result follows by combining  all the facts.
For $\alpha$ be the endpoints, due to the extremity property the proof is  easy.

For our case $\LL_f\subset \R^d$ is a compact convex set. If $\alpha\in {\rm ri}(\LL_f)$, it is possible to follow the line in \cite{JJOP2010} to give a lower bound. If $\alpha\in \LL_f$ is an extreme point,
the argument in \cite{JJOP2010} also works.
But now for $\alpha$ which is not in  the  relative interior and nor an extreme point, it seems quite hard to construct such an ergodic sequence $\mu_n$ which satisfies $\int A_n f d\mu_n=\alpha.$ The extremity argument is also not available in this case.

We will adapt a measure concatenation technique appeared in \cite{BQ} to construct a Moran subset  $M\subset X_\alpha$, on which we support a suitable measure, whose dimension can be estimated. By this way we obtain a unified way to treat the lower bound.
The construction of the Moran set also uses ideas from \cite{FFW}. See Section \ref{thm1-lower} for the technical  details.

 }
\end{remark}

\begin{remark}
{\rm
The second difficulty comes from the proof of the lower bound in Theorem \ref{main-2} when $\alpha \in A.$ For the one dimensional case  considered in \cite{JJOP2010}, $A$ is still a compact interval $A=[a_1,a_2].$
  They  deal with the lower bound by considering the boundary points and interior points separately.
When $\alpha$ is an interior point, by a suitable reduction, they can still use the argument presented in the former case. To deal with  the boundary point case  they adapt an approach appeared in \cite{GR2009}, which in essence is very closed to  the one used in \cite{BQ}.

Go back to higher dimensional case, we face with essentially the same difficulty, that is, for $\alpha\in A$ which is neither in the relative interior or not an extreme point of $A$, we can not apply the argument in \cite{JJOP2010} directly.

By carefully examining the proof in \cite{JJOP2010}, we find that indeed it can be modified a bit  to give a unified proof of the lower bound for all $\alpha\in A$.
}
\end{remark}

%\begin{remark}
%Our method is different from that in JJOP\cite{JJOP2010}. In JJOP\cite{2010}, they construct $n$-th Bernoulli measure whose dimension is well approximated to the lower dimension of $\alpha$-level set when $\alpha$ in the interior of  $\mathcal{L}_{f}$.  It will be almost trivial for the case that when $\alpha$ on the boundary of $\mathcal{L}_{f}$ in dimension $1$.
%
%In higher dimensions, the case for boundary will be much more complicated since it is very possible that different parts of the boundary vary in dimensions. Also it seems that there is not a direct proof to generalize JJOP's methods dealing with the interior case. We give a unified approach to estimate the lower bound for Hausdroff dimension of level set and our approach is based on to construct Moran set in the  $\alpha$ Level set.
%
%
%\end{remark}

%\begin{remark}
%As pointed in JJOP\cite{JJOP2010}, whether the assumption in {\bf Theorem 2} can be removed or not is still an open problem.  Urbanski \cite{Urbanski1996} proved it is true when the inverse branch of $T$ is $C^{1+\alpha} (0<\alpha<1)$ with some geometric conditions.  In some sense, this problem is closely related to whether the least upper bound of dimensions of hyperbolic measure coincides with the zeros of topological pressure. In this direction, for uniformly hyperbolic dynamic system Gatzouras and Peres \cite{GP1997} proved this is true even for $C^{1}$ case. However, it is still be open for the non-uniform hyperbolic case in $C^{1}$ condition.
%\end{remark}

The rest of this  note is organized as follows.  In Section \ref{preliminary} we give  some preliminary results  and lemmas which are needed for the proof. In Section \ref{thm1-lower} we prove the lower bound of Theorem \ref{main-1}. In Section \ref{thm1-upper} we prove the upper bound of Theorem \ref{main-1}.
In Section \ref{thm2} we prove Theorem \ref{main-2}. In Section \ref{coro-1} we prove Corollary \ref{main-3}.
% Finally in Section \ref{example} we give an example to demonstrate our main results.

%%%%%%%%%%%%%%%%%%%%%%%%
\section{Preliminaries}\label{preliminary}

In this section, we will give the notations and the lemmas needed in the proof.

Assume $T:X\to X$ is a topological dynamical system.
Denote by $\M(X,T)$ the set of all invariant probability measures and ${\mathcal E}(X,T)$ the set of all ergodic probability measures. Given $\mu\in\M(X,T)$, let $h(\mu,T)$ be the metric entropy of $\mu$.
Given  $f:X\to\R^d$ continuous, we write
$$
S_n f(x)=\sum_{j=0}^{n-1}f(T^jx) \ \ \text{ and }\ \  A_nf(x)=\frac{S_nf(x)}{n}.
$$
$S_nf(x)$ is called the {\it ergodic sum} of $f$ at $x$ and $A_nf(x)$ is called the {\it ergodic average} of $f$ at $x$.

Recall that  $\mathcal{A}=\{1,2\dots m\}$ and $\Sigma=\mathcal{A}^{\mathbb{N}}$.  Write
$
\Sigma_n=\{w=w_1\cdots w_n: w_i\in\mathcal{A}\}.
$
 For $\omega=\{\omega_{n}\}_{n=1}^{\infty}\in \Sigma$, write $\omega|_n=\omega_1\cdots \omega_n$.
 For $w\in\Sigma_n$  define the cylinder $[w]:=\{\omega: \omega|_n=w\}$.
  For a continuous potential $f:\Sigma\rightarrow \R^{d}$,  define the {\it $n$-th variation} of $f$  as
$$
\var_{n}f=\sup_{\omega|n=\tau|n}|f(\omega)-f(\tau)|,
$$
where $|\cdot|$  is the Euclidean norm in $\R^{d}$. Given $f: \Sigma\rightarrow \R^d$ continuous, let
$\|f\|:=\underset{\tau\in\Sigma}{\sup}|f(\tau)|$.  For $f: \Lambda\rightarrow \R^d$ continuous we define $\|f\|$ similarly.
 We have the following standard result:

\begin{lemma}\label{variation}
 If $f:\Sigma\to\R^d$ is continuous, then $\var_n A_nf \to 0.$
\end{lemma}

Consider the projection $\Pi: \Sigma\to \Lambda$. Let
$\tilde \Lambda:=\{x\in \Lambda: \# \{\Pi^{-1}(x)\}=2\}$.
In other words $\tilde \Lambda$ is the set of such $x$ with two codings. By our assumption on $I_j$, we know that both  $\tilde \Lambda$ and $\Pi^{-1}\tilde \Lambda$ are  at most countable. Moreover
\begin{equation}\label{coding}
\Pi^{-1}\tilde \Lambda\subset \{\omega: \omega= wm^\infty \text{ or } \tilde w 1^\infty\}.
\end{equation}

Then it is seen that
$$
\Pi: \Sigma\setminus \Pi^{-1}(\tilde \Lambda)\to \Lambda\setminus \tilde \Lambda
$$
is a bijection. We will need this fact in the proof of the lower bound of Theorem \ref{main-1}.

For $w=w_1\cdots w_n$, write $I_w=T_{w_1}\circ \dots \circ T_{w_n}[0,1].$ Especially for $\omega\in\Sigma$, we write $I_n(\omega)=I_{\omega|_n}$.
Let $D_{n}(\omega)=\diam(I_{n}(\omega))$.
Recall that we have defined
  $g(\omega):=-\log T'_{\omega_1}\Pi(\sigma \omega)$ and
$$
\tilde{\Sigma}=\{\omega\in\Sigma:\liminf\limits_{n\rightarrow \infty }A_ng(\omega)>0\}.
$$
$D_n(\omega)$ can be estimated via $A_ng(\omega)$ by the following lemma:
\begin{lemma}[\cite{Urbanski1996,JJOP2010}]\label{appro}
Under the assumption on $T$, $D_{n}(\omega)$  converges to $0$ uniformly. Moreover
$$
\lim\limits_{n\rightarrow \infty}\sup\limits_{\omega\in \Sigma}\left\{|-\frac{1}{n}\log D_{n}(\omega)-A_{n}g(\omega)|\right\}=0.
$$
\end{lemma}
 By this lemma we can understand that $\tilde \Sigma$ is the set  of such points $\omega$ such that the length of $I_n(\omega)$ tends to 0 exponentially. To simplify the notation we write $\tilde \lambda_n(\omega)=-\log D_n(\omega)/n.$

 Given $\mu\in\M(\Sigma,\sigma)$,  let
 $\lambda(\mu,\sigma):=\int g d\mu$ be the Lyapunov exponent of $\mu.$ Similarly given $\mu\in\M(\Lambda,T)$, let
 $\lambda(\mu,T):=\int \log|T^\prime| d\mu$ be the Lyapunov exponent of $\mu.$ For a $\mu\in\M(\Sigma,\sigma)$, we let $\Pi_*\mu=\mu\circ\Pi^{-1}$ in this paper.

The following lemma, which is a combination of Lemma 2 and Lemma 3 in \cite{JJOP2010}, is very useful in our proof.

\begin{lemma}\label{basic}
For any $\mu\in \mathcal{M}(\Sigma,\sigma)$, there exists  a sequence of ergodic measures $\{\mu_n:n\ge 1\}$ such that
  $\mu_n\rightarrow \mu$   in the weak star topology and
 $$
 h(\mu_n,\sigma)\rightarrow h(\mu,\sigma);
 \ \ \ \lambda(\mu_n,\sigma)\rightarrow \lambda(\mu,\sigma).
 $$
\end{lemma}

We remark that from their proof each  ergodic measure $\mu_n$ is continuous, i.e. $\mu_n$ has no atom.

\section{Lower bound for Theorem \ref{main-1}}\label{thm1-lower}

This section is devoted to the proof of the lower bound for Theorem \ref{main-1}. At first we show that the equality in Theorem \ref{main-1} makes sense.

\begin{lemma}
Assume that $X_{\alpha}\cap\tilde{\Sigma}\neq \emptyset$, then there exists a $\mu\in \mathcal{M}(\Sigma,\sigma)$ such that $\int fd\mu=\alpha$, and $\lambda(\mu,\sigma)>0$.
\end{lemma}
\begin{proof}
In fact, for any $\omega\in X_{\alpha}\cap\tilde{\Sigma}$, we can take a subsequence  $\{n_{k}\}_{k=1}^{\infty}$ such that
\[\liminf_{n\rightarrow\infty} A_{n}g(\omega)=\lim_{k\rightarrow \infty} A_{n_k}g(\omega)=\lim_{k\rightarrow\infty}\int g d\left(\frac{1}{n_k}\sum_{i=0}^{n_k-1}\delta_{\sigma^{i}\omega}\right).\]
Without loss of generality, we can assume that
$\frac{1}{n_k}\sum\limits_{i=0}^{n_k-1}\delta_{\sigma^{i}\omega}$ converges to some $\mu$ in $\mathcal{M}(\Sigma,\sigma)$ weakly. Since $\omega\in\tilde \Sigma$, we have $\lambda(\mu,\sigma)>0$. Since $\omega\in X_\alpha$ we have  $\int fd\mu=\alpha$.
\end{proof}

Recall that $\tilde \Lambda$ is the set of $x\in \Lambda$ which has two codings. The lower bound is a direct consequence of the following lemma:

\begin{lemma}{\label{key}}
Given $\mu \in \mathcal{M}(\Sigma,\sigma)$, such that $\int fd\mu=\alpha$ and $\lambda(\mu,\sigma)>0$. We can construct a Moran set $M\subset (X_{\alpha}\cap \tilde{\Sigma})\setminus \Pi^{-1}\tilde \Lambda$ together with a probability measure $\nu$ supported on it such that
$$
\underline{d}_{\Pi_{*}\nu}(x)\geq \frac{h(\mu,\sigma)}{\lambda(\mu,\sigma)}
$$
for all $x\in \Pi(M)$, where $\underline{d}_{\Pi_{*}\nu}(x)$ is the lower local dimension of $\Pi_*\nu$ at $x$.
\end{lemma}

\noindent {\bf Proof of Theorem \ref{main-1}: lower bound.}\ It follows from Lemma \ref{key} that for any $\mu \in \mathcal{M}(\Sigma,\sigma)$ such that $\int fd\mu=\alpha$ and $\lambda(\mu,\sigma)>0$ we can construct a measure $\nu$ such that $\nu(M)=1$ and satisfies
$$
\underline{d}_{\Pi_{*}\nu}(x)\geq \frac{h(\mu,\sigma)}{\lambda(\mu,\sigma)}
$$
for $\Pi_{*}\nu$ a.e. $x$. This implies that
$\dim_H \Pi_*\nu\ge {h(\mu,\sigma)}/{\lambda(\mu,\sigma)}$(see for example \cite{F}).
Since $\Pi:\Sigma\setminus \Pi^{-1}\tilde \Lambda\to \Lambda\setminus \tilde \Lambda$ is bijection and $M\cap \Pi^{-1}\tilde \Lambda=\emptyset$, we conclude that
$\Pi_\ast\nu(\Pi (M))=\nu(M)=1.$ Then
$$
\dim_{H}\Pi(X_{\alpha}\cap\tilde{\Sigma})\ge\dim_H \Pi(M)\geq d_{H}\Pi_{*}\nu \geq \frac{h(\mu,\sigma)}{\lambda(\mu,\sigma)}.
$$
Take a supremum we get the desired lower bound.
\hfill $\Box$

\medskip

It remains  to prove Lemma \ref{key}. As mentioned in the introduction, we will apply Lemma \ref{basic} to get some building blocks. Then we will concatenate them in such a way that we can construct a Moran set $M$ sitting inside the level set and supporting our limit measure. On the other hand during the concatenating process we can also control the size of the measure in cylinders, that finally we can also get the desired local dimension estimation.

\medskip

\noindent {\bf Proof of Lemma \ref{key}.}\
  By Lemma \ref{variation} and Lemma \ref{appro}, we can choose a decreasing sequence $\epsilon_i\downarrow 0$   such that for all $n\geq i$,
\begin{equation}\label{var}
\var_n A_n f<\epsilon_i,\ \ \ \var_n A_n g<\epsilon_i\ \ \text{ and }\ \  \ |\tilde{\lambda}_n(\omega)-A_n g(\omega)|<\epsilon_i (\forall \omega\in\Sigma).
\end{equation}
By Lemma \ref{basic} we can pick a sequence of  $\nu_i\in\mathcal{E}(\Sigma,\sigma)$, such that
\begin{equation}\label{control}
|\int f d \nu_i-\alpha|<\epsilon_i,\
|h(\nu_i,\sigma)-h(\mu,\sigma)|<\epsilon_i\ \text{ and }\
|\lambda(\nu_i,\sigma)-\lambda(\mu,\sigma)|<\epsilon_i.
\end{equation}
Since $\nu_i$ is ergodic,  for $\nu_i$ a.e. $\omega$,
\begin{equation}\label{block}
A_n f(\omega)\to \int f d \nu_i,  \
A_n g(\omega)\to \lambda(\nu_i,\sigma)  \text{ and }
-\frac{\log \nu_i[\omega|_n]}{n}\to h(\nu_i,\sigma).
\end{equation}

Fix $\delta>0$. Since $\nu_i$ is continuous as we remarked after Lemma \ref{basic}, there exists $\ell_i\ge i$ such that
$\nu_i(\bigcup_{j=1}^m [j^{\ell_i}])\le \delta/2.$
By  Egorov's theorem,  there exists $\Omega'(i)\subset \Sigma$ such that $\nu_i(\Omega'(i))>1-\delta/2$ and \eqref{block}
 holds uniformly on $\Omega'(i)$. Then there  exists   $l_i\geq \ell_i\ge i$ such that  for all  $n\geq l_i$  and $\omega\in\Omega'(i)$, we
have
\begin{equation}\label{estimation}
\begin{cases}
|A_n f(\omega)- \int f d \nu_i|<\epsilon_i \\
|A_n g(\omega)- \lambda(\nu_i,\sigma)|< \epsilon_i \\
|-{\log\nu_i[\omega|_n]}/{n} - h(\nu_i,\sigma)|< \epsilon_i
\end{cases}
\end{equation}
Let
$$
\Sigma(i)=\{\omega|_{l_{i}}\ |\ \omega\in \Omega'(i)\}\setminus\{1^{l_i},\cdots, m^{l_i}\}.
$$
Let ${\Omega}(i)=\bigcup_{w\in \Sigma(i)}[w]$. Then
$$
\nu_i({\Omega}(i))\ge \nu_i(\Omega^\prime(i))-\nu_i(\bigcup_{j=1}^m [j^{l_i}]) \ge 1-\delta/2 -\delta/2=1-\delta.
$$

It is seen that we can take $l_i$ such that  $l_i\uparrow \infty$ and still satisfies all the above property.
Let $N_i=l_{i+2}$, $i\geq 1$. Let
$$
M=\overset{\infty}{\underset{i=1}{\prod}}\overset{N_i}{\underset{j=1}{\prod}}\Sigma(i).
$$
By the definition of $\Sigma(i)$ and \eqref{coding}, it is ready to see that
$M\cap \Pi^{-1}\tilde \Lambda =\emptyset.$ In the following we will construct a measure supporting on it and show that $M\subset X_\alpha\cap \tilde \Sigma.$

Relabel the following sequence
$$
\underbrace{l_1\cdots l_1,}_{N_1}\cdots,\underbrace{l_i\cdots l_i,}_{N_i}\cdots
$$
 as $\{l^*_i:i\ge 1\}$.
Relabel the following sequence
$$
\underbrace{\Sigma(1)\cdots \Sigma(1),}_{N_1}\cdots,\underbrace{\Sigma(i)\cdots
\Sigma(i),}_{N_i}\cdots
$$
 as $\{\Sigma^{*}(i): i\ge 1\}$.
Accordingly  we get $\{{\Omega}'^*(i)\}$, $\{{\Omega}^*(i)\}$, $\{\nu^*_i)\}$, $\{\epsilon^*_i\}$. Let
$n_k=\underset{i=1}{\overset{k}{\sum}}l_i^*$. For any $n>0$, there  exists $J(n)\in\N$ such that
$\underset{i=1}{\overset{J(n)}{\sum}}l_i^*\leq n<\underset{i=1}{\overset{J(n)+1}{\sum}}l_i^*$.
There also exists $r(n)\in\N$ such that $\underset{i=1}{\overset{r(n)}{\sum}}N_i\leq
J(n)<\underset{i=1}{\overset{r(n)+1}{\sum}}N_i$. It is seen that
\begin{equation}\label{J-n}
J(n)\leq J(n+1)\leq J(n)+1,\ l_{J(n)+1}^*=l_{r(n)+1}\ \text{ and }\ l_{J(n)+2}^*\leq l_{r(n)+2},
\end{equation}
then, for $j=1,2$, $$\frac{l_{J(n)+j}^*}{\underset{i=1}{\overset{J(n)}{\sum}}l_i^*}\leq
\frac{l_{r(n)+j}}{N_{r(n)}l_{r(n)}}=\frac{l_{r(n)+j}}{l_{r(n)+2}l_{r(n)}}.$$
Since $l_i$ is increasing to $\infty$,
we also have
\begin{equation}\label{l-i-basic}
{\underset{i=1}{\overset{J(n)+1}{\sum}}l_i^*}/{\underset{i=1}{\overset{J(n)}{\sum}}l_i^*}\to 1\ \ \text{ and }\ \ \ {l^*_{J(n)+j}}/{\underset{i=1}{\overset{J(n)}{\sum}}l_i^*}\to 0, j=1,2.
\end{equation}

At first  we define a probability $\nu$ supported on $M$. For each  $w\in \Sigma^*(i)$ define
$$
\rho^i_w=\frac{\nu_i^\ast[w]}{\nu_i^\ast(\Omega^*(i))}.
$$
It is seen that $\sum_{w\in \Sigma^*(i)}\rho^i_w=1.$
Write
$
{\mathcal C}_n:=\{[w]:w\in \prod_{i=1}^n \Sigma^\ast(i)\}.
$
It is seen that $\sigma({\mathcal C}_n: n\ge 1)$ gives the Borel-$\sigma$ algebra in $M.$
For each $w=w_1\cdots w_n\in {\mathcal C}_n$ define
$$
\nu([w])=\prod_{i=1}^n \rho^i_{w_i}.
$$
Let $\nu$ be the Kolmogorov extension of $\nu$ to all the Borel sets.   By the construction it is seen that $\nu$ is supported on $M.$

Next   we show that $M\subset X_\alpha \cap \tilde\Sigma$.
Write $n_0=0$ and  $n_i=\sum_{j=1}^i l_i^\ast$ for $i\ge 1$. Fix any $\omega\in M$. By the construction we have $\sigma^{n_{i-1}}\omega \in [w]$ for some $w\in\Sigma^*(i)$, consequently there exists $\omega^i\in \Omega^{\prime \ast}(i)\cap [w]$ such that \eqref{estimation} holds. So we have
\begin{align*}
&|S_n f(\omega)-n\alpha|\\
=&|\underset{i=1}{\overset{J(n)}{\sum}}l_i^*(A_{l_i^*}f
(\sigma^{n_{i-1}}\omega)-\alpha)
+S_{n-n_{J(n)}}f (\sigma^{n_{J(n)}}\omega)-(n-n_{J(n)})\alpha| \\
\leq&
|\underset{i=1}{\overset{J(n)}{\sum}}l_i^*\left(A_{l_i^*}f (\sigma^{n_{i-1}}\omega)-A_{l_i^*}f (\omega^i)+
A_{l_i^*}f (\omega^i)-\int f d \nu_i^*+\int f d \nu_i^*-\alpha\right)|\\
&+l_{J(n)+1}^*(\alpha+||f||)\\
\le&\sum_{i=1}^{J(n)} l_i^*(\epsilon_i^\ast+\epsilon_i^\ast+\epsilon_i^\ast)+l_{J(n)+1}^*(\alpha+||f||),
\end{align*}
where for the last inequality we use \eqref{var}, \eqref{control} and \eqref{estimation}.
  Then
  $$
  \frac{|S_n f(\omega)-n\alpha|}{n}\leq
\frac{3\underset{i=1}{\overset{J(n)}{\sum}}l_i^*\epsilon_i^*}{\underset{i=1}{\overset{J(n)}{\sum}}l_i^*}+
\frac{l_{J(n)+1}^*(\alpha+||f||)}{\underset{i=1}{\overset{J(n)}{\sum}}l_i^*}.
$$
By  \eqref{l-i-basic} and the fact that $\epsilon_i^*\downarrow 0$ we conclude that
$A_n f(\omega)\rightarrow\alpha$, which implies $x\in X_\alpha.$ Thus $M\subset X_\alpha.$

Now we check that $\underset{n\rightarrow\infty}{\liminf}A_ng(\omega)>0$ for all $\omega\in M$. Let $\omega^i$ defined as above, then
\begin{align*}
A_ng(\omega)&=\frac{1}{n}\underset{i=1}{\overset{J(n)}{\sum}}l_i^*A_{l_i^*}g (\sigma^{n_{i-1}}\omega)+
\frac{1}{n}S_{n-n_{J(n)}}g (\sigma^{n_{J(n)}}\omega)\\
&\geq
\frac{1}{n}\underset{i=1}{\overset{J(n)}{\sum}}l_i^*\huge[A_{l_i^*}g (\sigma^{n_{i-1}}\omega)-A_{l_i^*}g
(\omega^i)+A_{l_i^*}g
(\omega^i)-\lambda(\nu_i^*,\sigma)\\
&+\lambda(\nu_i^*,\sigma)-\lambda(\mu,\sigma)+\lambda(\mu,\sigma)\huge]-
\frac{1}{n}l_{j(n)+1}^*||g|| \\
&\geq
\frac{1}{n}\underset{i=1}{\overset{J(n)}{\sum}}l_i^*(\lambda(\mu,\sigma)-3\epsilon_i^*)-
\frac{1}{n}l_{j(n)+1}^*||g|| \\
&\geq
\frac{\underset{i=1}{\overset{J(n)}{\sum}}l_i^*(\lambda(\mu,\sigma)-3\epsilon_i^*)}
{\underset{i=1}{\overset{J(n)+1}{\sum}}l_i^*}-\frac{l_{j(n)+1}^*||g||}{\underset{i=1}{\overset{J(n)}{\sum}}l_i^*},
\end{align*}
where for the second inequality we again use \eqref{var}, \eqref{control} and \eqref{estimation}.  Now by \eqref{l-i-basic}
we get $\underset{n\rightarrow\infty}{\liminf}A_ng(\omega)\geq \lambda(\mu,\sigma)>0$. Thus $\omega\in \tilde \Sigma$ and we conclude that $M\subset \tilde \Sigma.$

Finally we  compute the local dimension of $\Pi_\ast\nu$.
we will  show that  for all $x\in \Pi(M)$
$$
\underset{r\downarrow0}{\liminf}\frac{\log \Pi_*\nu (B(x,r))}{\log r}\geq
 \frac{h(\mu,\sigma)}{\lambda(\mu,\sigma)}.
 $$

Fix  $\omega\in M$. At first we find a lower bound for $D_n(\omega)$. Define $n_i$ and $\omega^i$ as before.
Recall  that $D_n(\omega)=e^{-n\tilde \lambda_n(\omega)}$. By \eqref{var} we have
\begin{align*}
&  n\tilde \lambda_n(\omega)\\
\leq& n(A_ng(\omega)+\epsilon_{J(n)}^\ast)\\
\leq&
 \underset{i=1}{\overset{J(n)}{\sum}}l_i^*(A_{l_i^*}g (\sigma^{n_{i-1}}\omega)+\epsilon_i^*)
 +(n-n_{J(n)})\left(A_{n-n_{J(n)}}g (\sigma^{n_{J(n)}}\omega)+\epsilon_{J(n)}^*\right)\\
 \leq&
\underset{i=1}{\overset{J(n)}{\sum}}l_i^*\big\{A_{l_i^*}g(\sigma^{n_{i-1}}\omega)-A_{l_i^*}g(\omega^i)+A_{l_i^*}g
(\omega^i)-\lambda(\nu_i,\sigma)+\\
&\lambda(\nu_i,\sigma)-\lambda(\mu,\sigma)
+\lambda(\mu,\sigma)+\epsilon_i^*\big\}
+l_{J(n)+1}^*(||g||+\epsilon_{J(n)}^*)\\
 \leq &\underset{i=1}{\overset{J(n)}{\sum}}l_i^*(\lambda(\mu,\sigma)+
4\epsilon_i^*)+
l_{J(n)+1}^*(||g||+\epsilon_{J(n)}^*)=:\rho(n).
\end{align*}
Then
$D_n(\omega)\geq e^{-\rho(n)}.$ It is seen that $\rho(n)$ is increasing.

Now fix $x\in \Pi(M)$ and some $r>0$ small. Then there exists a unique $n=n_r$ such that
\begin{equation}\label{r}
e^{-\rho(n+1)}\le r< e^{-\rho(n)}.
\end{equation}
Consider the set of $n$-cylinders
$$
{\mathcal C}:=\{I_n(\omega): \omega\in M \text{ and } I_n(\omega)\cap B(x,r)\ne \emptyset\}.
$$
By the bound $D_n(\omega)\ge e^{-\rho(n)}$, the above  set consists of at most three cylinders, i.e. $\#{\mathcal C}\le 3$.

Choose $\omega\in M$ such that  $I_n(\omega)\in {\mathcal C}$. Write $\omega|_n=w_1\cdots w_{J(n)} v$, then $w_i\in \Sigma^*(i)$ and $v$ is a prefix of some $\tilde v\in \Sigma^*(J(n)+1)$.
Then
\begin{align*}
 \Pi_*\nu(I_n(\omega))=\nu[\omega|_n]
=&\underset{i=1}{\overset{J(n)}{\prod}}\frac{\nu_i^*[w_i]}{\nu_i^*(\Omega^*(i))}\cdot\frac{\nu_{J(n)+1}^*[v]}{\nu_{J(n)+1}^*(\Omega^*(J(n)+1))}\\
& \leq
(1-\delta)^{-J(n)-1}\underset{i=1}{\overset{J(n)}{\prod}}\nu_i^*[w_i].
\end{align*}
Then we conclude that $\Pi_*\nu(B(x,r))\le 3(1-\delta)^{-J(n)-1}\underset{i=1}{\overset{J(n)}{\prod}}\nu_i^*[w_i].$ Consequently
\begin{align*}
&\quad \log \Pi_*\nu(B(x,r))\\
&\leq
-\underset{i=1}{\overset{J(n)}{\sum}}l_i^*\left(-\frac{\log\nu_i^*[w_i]}{l_i^*} \right)
-(J(n)+1)\log(1-\delta)+\log 3\\
&\leq
-\underset{i=1}{\overset{J(n)}{\sum}}l_i^*(h(\mu,\sigma)-2\epsilon_i^*)-(J(n)+1)\log(1-\delta)+\log 3,
\end{align*}
where for the second inequality  we use \eqref{control} and \eqref{estimation}.
Notice  that $r\rightarrow 0$ if and only if $n\rightarrow\infty$. By \eqref{J-n} we have $J(n+1)\le J(n)+1.$ Together with \eqref{r} and \eqref{l-i-basic} we get
\begin{align*}
& \quad \underset{r\downarrow 0}{\liminf}\frac{\log\Pi_* \nu(B(x,r))}{\log r} \\
& \geq
\underset{n\rightarrow\infty}{\lim}
\frac{\underset{i=1}{\overset{J(n)}{\sum}}l_i^*(h(\mu,\sigma)-2\epsilon_i^*)+(J(n)+1)\log(1-\delta)-\log 3}
{\underset{i=1}{\overset{J(n+1)}{\sum}}l_i^*(\lambda(\mu,\sigma)+4\epsilon_i^*)
+l_{J(n+1)+1}^*(||g||+\epsilon_{J(n+1)}^*)}\\
&=\frac{h(\mu,\sigma)}{\lambda(\mu,\sigma)}.
\end{align*}
Then the result follows.\hfill $\Box$

%%%%%%%%%%%%%%%%%%%%%%%%%%%

\section{Upper bound for Theorem \ref{main-1}}\label{thm1-upper}
  The proof of the upper bound    is essentially the same with that given in \cite{JJOP2010}. We  include it for completeness.

  For each  $n\in\N$ define
  $$
\tilde \Sigma(n)=\{\omega\in\Sigma: \liminf_{n\to\infty }A_n g(\omega)\ge 1/n\}.
$$
 For $\alpha\in \LL_f, N, k,n\in \N$ define
$$
X(\alpha, N,k,n)=\{\omega\in \tilde{\Sigma}(n): A_{p}f(\omega)\in B(\alpha,1/k), \textit{for all } p\geq N\}.
$$
We need the following lemma:
\begin{lemma}[\cite{JJOP2010}]\label{JJOP}
For $ \epsilon>0$ sufficiently small, $k,n\in \N$ sufficiently big and $N\in\N$, we can find a measure $\mu \in \mathcal{M}(\Sigma,\sigma)$, such that $\int f d\mu \in B(\alpha,2/k)$, $\lambda(\mu,\sigma)>1/n$ and
\[
\dim_{H}\Pi(X(\alpha, N,k,n))\leq \frac{h(\mu,\sigma)}{\lambda(\mu,\sigma)}+\epsilon.
\]
\end{lemma}

At first we have
$$
X_{\alpha}\cap \tilde{\Sigma}=\bigcup_{n=1}^{\infty}X_{\alpha}\cap \tilde{\Sigma}(n).
$$

Fix $n\in \N.$
For all $k\in\N$  we  also have
$$
X_{\alpha}\cap \tilde{\Sigma}(n)\subset \bigcup_{N\ge k} X(\alpha, N,k,n).
$$
Consequently  for all $k\in\N$
\begin{eqnarray*}
 \dim_H\Pi(X_{\alpha}\cap \tilde{\Sigma}(n))
\le \sup_{N\ge k} \dim_H \Pi(X(\alpha, N,k,n))=:\beta(k).
\end{eqnarray*}

 For each $k$, choose $X(\alpha,N_k,k,n)$ such that
 $$
 \beta(k)\le \dim_H\Pi(X(\alpha,N_k,k,n))+\epsilon/2.
 $$
 By Lemma \ref{JJOP}, we can choose $\mu_k\in\M(\Sigma,\sigma)$ such that $\int f d\mu_k\in B(\alpha,2/k)$ $\lambda(\mu_k,\sigma)\ge 1/n$ such that
 \[
\dim_{H}\Pi(X(\alpha, N_k,k,n))\leq \frac{h(\mu_k,\sigma)}{\lambda(\mu_k,\sigma)}+\epsilon/2.
\]
Let $\mu_*$ be any weak star limit of $\mu_k$, it is clear that $\int f d\mu_*=\alpha$ and $\lambda(\mu_*,\sigma)\ge 1/n.$ Without loss of generality we assume $\mu_k\to \mu_*$ weakly. We note that in the symbolic case the metric entropy is upper semi-continuous, then we have
\begin{eqnarray*}
&& \dim_H\Pi(X_{\alpha}\cap \tilde{\Sigma}(n))
\le \limsup_k\frac{h(\mu_k,\sigma)}{\lambda(\mu_k,\sigma)}+\epsilon\\
&\le& \frac{ h(\mu_\ast,\sigma)}{\lambda(\mu_\ast,\sigma)}+\epsilon\\
&\le& \sup_{\mu\in\M(\Sigma,\sigma)}\left\{\frac{ h(\mu,\sigma)}{\lambda(\mu,\sigma)}: \int fd\mu=\alpha; \lambda(\mu,\sigma)>0\right\}+\epsilon.
\end{eqnarray*}
Since  we have
$$
\dim_H \Pi(X_\alpha\cap \tilde \Sigma)=\sup_{n}\dim_H\Pi(X_{\alpha}\cap \tilde{\Sigma}(n)) $$
and  $\epsilon>0$ is arbitrary, the desired upper bound holds.
\hfill $\Box$

\section{Proof for Theorem 2}\label{thm2}

We want to apply Theorem \ref{main-1} to this situation. At first we define $f=F\circ \Pi$, then $f$ is a continuous function on $\Sigma,$ moreover in this case $
\Lambda_\alpha=\Pi(X_\alpha).
$
Assume $I=\{x_{i_1},\cdots,x_{i_k}\}$ be the set of parabolic fixed points. Then $x_{i_l}$ has coding $\overline{i_l}:=i_l^\infty.$ So
$$
A={\rm Co}\{F(x_{i_l}):l=1,\cdots,k\}={\rm Co}\{f(\overline{i_l}):l=1,\cdots,k\}.
$$

To apply Theorem \ref{main-1} we need to know what is the relation between $X_\alpha$ and $\tilde \Sigma.$ The following lemma proved in \cite{JJOP2010} make this relation clear:
\begin{lemma}[\cite{JJOP2010}]\label{b}
Let $\{n_{j}\}$ be a subsequence of $\N$. If $\lim_{j\rightarrow \infty}A_{n_j}g(\omega)=0$ for some $\omega \in \Sigma$,  then we have that $\lim_{n \rightarrow \infty}A_{n_j}f(\omega)\in A$, if the limit exists. In particular, this shows that $\liminf_{n \rightarrow \infty}A_{n}g(\omega)=0$ means that $\lim_{n\rightarrow \infty}A_{n}f(\omega)\in A$, if the limit exists.
\end{lemma}
We remark that although in present  case $f$ is a vector valued function, the proof is essentially the same, so we omit it.

We also need the following property which establishes the relation between  $\M(\Sigma,\sigma)$ and $\M(\Lambda,T).$

\begin{lemma}[\cite{FLW2002}]\label{A}
Let $X_i$, $i=1,2$ be compact metric spaces and let $T_i: X_i\rightarrow X_i$ be continuous.
 Suppose $\Pi:X_1\rightarrow X_2 $ is a continuous surjection such
that the following diagram commutes:
\[\begin{array}{cccc}
X_1 &
\stackrel{T_1}{\longrightarrow} &
X_1  \\
\Big\downarrow \vcenter{%
\rlap{$\scriptstyle{\Pi}$}}& &
\Big\downarrow\vcenter{%
\rlap{$\scriptstyle{\Pi}$}}\\
X_2 &\stackrel{T_2} {\longrightarrow}&
X_2.
\end{array}\]
Then  $\Pi_*: \mathcal{M}(X_1,T_1)\rightarrow \mathcal{M}(X_2,T_2)$  is surjective. If
furthermore, there is an integer $m>0$ such that $\Pi^{-1}(y)$  has at most $m$ elements
for each $y\in X_2$; then $$h(\mu,T_1)=h(\Pi_*\mu,T_2)$$
for each $\mu\in\mathcal{M}(X_1,T_1)$.
\end{lemma}

%Applying Lemma \ref{A} to $F$, we have that $$\underset{\mu\in \mathcal{M}(\Sigma,\sigma)}{\sup}\left\{\frac{h(\mu,\sigma)}{\lambda(\mu,\sigma)}\left|\int F\circ\Pi d \mu=\alpha\right. \right\}=\underset{\nu\in \mathcal{M}(\Lambda,T)}{\sup}\left\{\frac{h(\nu,T)}{\lambda(\nu,T)}\left|\int F d \nu=\alpha\right. \right\}.$$

We also remark that,  applying to our case,  it is easy to check that for any $\mu\in \mathcal{M}(\Sigma,\sigma)$ we have $\lambda(\Pi_*\mu,T)=\lambda(\mu,\sigma)$
\medskip

\noindent{\bf Proof of Theorem \ref{main-2} when $\alpha\not\in A.$}\  If $\alpha\not\in A$ and $\omega\in X_\alpha$, then $A_n f(\omega)\to \alpha$. By Lemma \ref{b}, we have
$\liminf_{n\to\infty}A_ng(\omega)>0.$ Thus $X_\alpha \subset \tilde \Sigma.$ Consequently
$$
\Lambda_\alpha=\Pi(X_\alpha)=\Pi(X_\alpha\cap \tilde\Sigma).
$$

Given $\mu\in\Sigma(\Lambda,T)$ such that $\int F d\mu=\alpha.$ We claim that $\lambda(\mu,T)=\int \log T^\prime  d\mu >0.$ In fact if otherwise, we have $\int \log T^\prime d\mu=0.$ However by our assumption $\log T^\prime(x)\ge 0$ and  $\log T^\prime(x)=0$ if and only if $x\in I.$ Thus $\mu$ is supported on $I$. But this obviously implies that
$$
\alpha=\int F d\mu=\sum_{x\in I} F(x)\mu(x)\in A,
$$
which is a contradiction.

Now by Lemma \ref{A}, it is ready to see that if $\alpha\not\in A$, then  Theorem \ref{main-2} is a direct consequence of Theroem \ref{main-1}.
\hfill $\Box$

 So it remains to prove Theorem \ref{main-2} when $\alpha\in A.$ The only nontrivial part is  the lower bound. As we will see soon, our solution is quite similar with that has given in the proof of Theorem \ref{main-1},
 we still need to concatenate measures and construct Moran sets, except in this case we will concatenate some narrow cylinders for certain stage.

 \medskip

\noindent{\bf Proof of Theorem \ref{main-2} when $\alpha\in A.$}\   The upper bound is trivial in this case. So we only prove the lower bound.
  By the assumption of Theorem \ref{main-2}, for any $\epsilon>0$, there  exists $\tilde\nu\in \M(\Lambda,T)$ with $\lambda(\tilde\nu,T)>0 $ such that $\frac{h(\tilde\nu,T)}{\lambda(\tilde\nu,T)}>\dim_{\text{H}}\Lambda-\epsilon.$ By Lemma 8 and the remark after it,
   there exists $\nu\in\mathcal{M}(\Sigma,\sigma)$ such that $\Pi_\ast(\nu)=\tilde\nu$ with $\lambda(\nu,\sigma)>0$ and
$\frac{h(\nu,\sigma)}{\lambda(\nu,\sigma)}>\dim_{\text{H}}\Lambda-\epsilon$.
By Lemma \ref{basic}, there exists a measure $\mu\in \mathcal{E}(\Sigma,\sigma)$ such that  $\lambda(\mu,\sigma)>0$ and $\frac{h(\mu,\sigma)}{\lambda(\mu,\sigma)}>\dim_{\text{H}}\Lambda-2\epsilon$.

We will construct a Moran set
$M\subset X_{\alpha}$ and a probability measure $\eta$ such that $\eta(M)=1$ and $\dim_{\text{H}}\Pi_*\eta\geq \frac{h(\mu,\sigma)}{\lambda(\mu,\sigma)}$. This will end the proof since then $\Pi(M)\subset \Pi(X_\alpha)=\Lambda_\alpha$ and consequently  $\Pi_*\eta(\Lambda_\alpha)\ge \Pi_*\eta(\Pi(M))\ge \eta(M)=1.$ Thus
$$
\dim_H\Lambda_\alpha\ge \dim_{\text{H}}\Pi_*\eta\geq \frac{h(\mu,\sigma)}{\lambda(\mu,\sigma)}>\dim_{\text{H}}\Lambda-2\epsilon.
$$
Since $\epsilon$ can be arbitrarily small, we get the lower bound.

Now we begin to construct $M$ and $\eta.$
  Note that $A$ is a convex polyhedron. Since $\alpha\in A$, we can find extreme points of $A$ such that $\alpha$ is a strict convex combination of them. After relabeling the extreme points,  without loss of generality  we assume that  there exist real numbers $\{r_i\}_{i=1}^s$ such that
  \begin{equation}\label{convex}
  r_i> 0,\ \ \ \underset{i=1}{\overset{s}{\sum}}r_i=1 \ \ \text{ and  }\ \  \alpha=\underset{i=1}{\overset{s}{\sum}}r_i f(\overline{i}).
  \end{equation}

By Lemma \ref{variation} and Lemma \ref{appro}, we can choose a decreasing sequence $\epsilon_i^\prime\downarrow 0$   such that for all $n\geq i$,
\begin{equation}\label{var-2}
\var_n A_n f<\epsilon_i^\prime,\ \ \ \var_n A_n g<\epsilon_i^\prime\ \ \text{ and }\ \  \ |\tilde{\lambda}_n(\omega)-A_n g(\omega)|<\epsilon_i^\prime (\forall \omega\in\Sigma).
\end{equation}

 Assume $\int f d\mu=\beta.$
Since $\mu$ is ergodic,  for $\mu$ a.e. $\omega$,
$$
A_n f(\omega)\to \beta,  \
A_n g(\omega)\to \lambda(\mu,\sigma)  \text{ and }
-\frac{\log \mu[\omega|_n]}{n}\to h(\mu,\sigma).
$$

Fix $\delta>0$.
By  Egorov's theorem, we can choose another decreasing sequence $\tilde \epsilon_i\downarrow 0$ such that  there exists $\Omega'(i)\subset \Sigma$ with  $\mu(\Omega'(i))>1-\delta$ and    for all  $n\geq i$  and $\omega\in\Omega'(i)$, we
have
\begin{equation}\label{estimation-1}
\begin{cases}
|A_n f(\omega)- \beta|<\tilde \epsilon_i \\
|A_n g(\omega)- \lambda(\mu,\sigma)|< \tilde \epsilon_i \\
|-{\log\mu[\omega|_n]}/{n} - h(\mu,\sigma)|< \tilde \epsilon_i
\end{cases}
\end{equation}
Let
$$
\Sigma(i)=\{\omega|_{{i}}\ |\ \omega\in \Omega'(i)\}.
$$
Let ${\Omega}(i)=\bigcup_{w\in \Sigma(i)}[w]$. Then
$$
\mu({\Omega}(i))\ge \mu(\Omega^\prime(i))\ge1-\delta.
$$

 Let $\epsilon_i=\max\{\epsilon_i^\prime, \tilde \epsilon_i\}.$
   We can select integer sequence $(k_i)_{i\geq1}$ such that
 \begin{equation}\label{k-i}
k_i\rightarrow \infty,\ \ k_i\epsilon_i\rightarrow 0,\ \ k_i\underset{1\leq j\leq s}{\min}\{r_j\}\geq1\ \text{ and }\ \underset{i}{\lim}\frac{k_i}{k_{i+1}}=1.
\end{equation}
For each $i\ge 1$ define
$$
v_j^i=\underbrace{j\cdots j}_{[r_jik_i]},\ \  (j=1\cdots s-1)\ \ \text{ and }\ \   v_s^i=\underbrace{s\cdots s}_{m_i}
$$ where  $[r_jik_i]$ is the integral part of  $r_jik_i$ and $m_i=ik_i-\underset{j=1}{\overset{s-1}{\sum}}[r_jik_i]$. Write $v^i= v_1^i\cdots v_s^i.$ Now we can define $M$ as
$$
M=\prod_{i=1}^\infty \left(\Sigma(i)\cdot v^i\right).
$$
Now  we define the  probability $\eta$ supported on $M$.  For each  $w\in \Sigma(i)$ define
$$
\rho^i_w=\frac{\mu([w])}{\mu(\Omega(i))}.
$$
It is seen that $\sum_{w\in \Sigma(i)}\rho^i_w=1.$
Write
$
{\mathcal C}_n:=\{[w]:w\in \prod_{i=1}^n \Sigma(i)\cdot v^i\}.
$
It is seen that $\sigma({\mathcal C}_n: n\ge 1)$ gives the Borel-$\sigma$ algebra in $M.$
For each $w=w_1v^1\cdots w_nv^n\in {\mathcal C}_n$ define
$$
\hat{\eta}([w])=\prod_{i=1}^n \rho^i_{w_i}.
$$
Let $\eta$ be the Kolmogorov extension of $\hat{\eta}$ to all the Borel sets.   By the construction it is seen that $\eta$ is supported on $M.$
In the following we only need to show that
\[
M\subset X_\alpha\ \ \ \text{
 and }\ \ \ \dim_H \Pi_*\eta\geq \frac{h(\mu,\sigma)}{\lambda(\mu,\sigma)}.
\]

  Let $n_0=0$ and
 $n_q=\underset{i=1}{\overset{q}{\sum}}i(1+k_i)$ for $q\ge 1$.
 Let $n_{q,0}=n_{q-1}+q$ and $n_{q,i}=n_{q-1}+q+\underset{l=1}{\overset{i}{\sum}}[r_lqk_q]$ for $i=1,\cdots,s-1$.

   Fix $\omega\in M$.   Note that since  $f$ is bounded and $\underset{q}{\lim}\frac{n_q-n_{q-1}}{n_q}=0$,  to show $\omega\in X_\alpha,$ we only need to consider the limit along the subsequence of $n_q$. By the construction we have $\sigma^{n_{i-1}}\omega \in [w]$ for some $w\in\Sigma(i)$, consequently there exists $\omega^i\in \Omega^{\prime}(i)\cap [w]$ such that \eqref{estimation-1} holds.
Write $t_{il}:=[r_lik_i]$ for $l=1,\cdots,s-1$ and $t_{is}=ik_i-(t_{i1}+\cdots+t_{i(s-1)}).$ We have
\begin{align*}
& \left|S_{n_q} f(\omega)-n_q \alpha\right|\\
=&\left|\sum_{i=1}^q \left(S_{n_i-n_{i-1}} f(\sigma^{n_{i-1}}\omega)-(n_i-n_{i-1})\alpha\right)\right|\\
=&\left|\sum_{i=1}^q \left[\left(S_{i} f(\sigma^{n_{i-1}}\omega)-i\alpha\right)+ \sum_{l=1}^s \left(S_{t_{il}} f(\sigma^{n_{i,l-1}}\omega)-t_{il}\alpha\right)\right]\right|\\
\le& \left|\underset{i=1}{\overset{q}{\sum}}i\left(A_i f (\sigma^{n_{i-1}}\omega)-A_i f(\omega^i)+A_i f(\omega^i)-\beta+\beta-\alpha\right)\right|+\\
&\sum\limits_{i=1}^{q}\left|\underset{l=1}{\overset{s}{\sum}}t_{il}\left(A_{t_{il}} f  (\sigma^{n_{i,l-1}}\omega)-A_{t_{il}} f(\overline{l})+f(\overline{l})-\alpha \right)\right|\ \ (A_{t_{il}} f(\overline{l})=f(\overline{l}))\\
\le& \underset{i=1}{\overset{q}{\sum}}i(2\epsilon_i+|\beta-\alpha|)+
\sum_{i=1}^q\underset{l=1}{\overset{s}{\sum}}t_{il}\epsilon_i+
\sum_{i=1}^q|\sum_{l=1}^s t_{il}f(\overline{l})- ik_i \alpha| \\
\le& \underset{i=1}{\overset{q}{\sum}}i(2\epsilon_i+|\beta-\alpha|)+\sum_{i=1}^q(ik_i\epsilon_i+
2(s-1)||f||) \\
=&\underset{i=1}{\overset{q}{\sum}}\left[i(2\epsilon_i+|\beta-\alpha|+k_i\epsilon_i)+2(s-1)||f||\right]
\end{align*}
where for the second inequality we use \eqref{var-2} and \eqref{estimation-1}, for the last inequality we use \eqref{convex}.
Now by \eqref{k-i} we have
\[\frac{1}{n_q}|S_{n_q} f(\omega)-n_q \alpha|\leq \frac{\underset{i=1}{\overset{q}{\sum}}[i(2\epsilon_i+|\beta-\alpha|+k_i\epsilon_i)+2(s-1)||f||]}{\sum_{i=1}^q ik_i}\to 0 \ \ \ (q\to\infty).
\]
Then we have
$\underset{q}{\lim} A_{n_q} f (\omega)=\alpha.$ This in turn implies $\lim_{n\to\infty}A_nf(\omega)=\alpha$. Thus $\omega\in X_\alpha.$ Consequently $M\subset X_\alpha.$

Now we will estimate the local dimension of $\Pi_*\eta.$ At first fix any $\omega\in M$, we estimate the length and $\eta$-measure of $I_{n_q}(\omega).$ We have
\begin{align*}
& n_{q}\tilde{\lambda}_{n_q}(\omega)\\
\le& S_{n_q}g(\omega)+n_q\epsilon_q\\
=&\sum_{i=1}^q S_{n_i-n_{i-1}} g(\sigma^{n_{i-1}}\omega)+n_q\epsilon_q\\
\le& \underset{i=1}{\overset{q}{\sum}}i\left(A_ig(\sigma^{n_{i-1}}\omega))-A_ig(\omega^i)+
A_ig(\omega^i)-\lambda(\mu,\sigma)+\lambda(\mu,\sigma)\right)+\\
&\sum_{i=1}^q\underset{l=1}{\overset{s}{\sum}}t_{il}(A_{t_{il}} g  (\sigma^{n_{i,l-1}}\omega)-A_{t_{il}} g(\overline{l}))+n_q\epsilon_q \ \ \ (g(\overline{l}))= A_{t_{il}}g(\overline{l})) =0)\\
\le& \underset{i=1}{\overset{q}{\sum}}i(\lambda(\mu,\sigma)+2\epsilon_i)+
\sum_{i=1}^qik_i\epsilon_i+\sum_{i=1}^qi(1+k_i)\epsilon_i\\
\le& \frac{q(q+1)}{2}\lambda(\mu,\sigma) +
5\sum_{i=1}^qik_i\epsilon_i =:\rho(q).
\end{align*}
Then
$D_{n_q}(\omega)\geq e^{-\rho(q)}.
$ Write $\omega|_{n_q}=w_1v^1\cdots w_qv^q$,
we also have
\begin{align*}
\quad \Pi_*\eta(I_{n_q}(\omega))= \eta([\omega|_{n_q}])
=\underset{i=1}{\overset{q}{\prod}}\frac{\mu[w_i]}{\mu(
{\Omega}(i))}
\leq (1-\delta)^{-q}\prod_{i=1}^q\mu([w_i]).
\end{align*}

Now fix any $x\in \Pi(M) $ and any small
 $r>0$. Then there exists a unique $q$ such that
 \begin{equation}\label{r-1}
 e^{-\rho(q+1)}\le r< e^{-\rho(q)}.
 \end{equation}

It is seen that  $B(x,r)$ can intersect at most 3  such $n_q$-level cylinders. Thus by \eqref{estimation-1}
\begin{align*}
 \log\Pi_*\eta(B(x,r))\leq&\log 3 -q\log(1-\delta)+
\underset{i=1}{\overset{q}{\sum}} \log\mu([w_i])\\
&\leq \log 3
-q\log(1-\delta)+\underset{i=1}{\overset{q}{\sum}}i\epsilon_i-\frac{q(q+1)}{2}h(\mu,\sigma).
\end{align*}

Then by \eqref{k-i} and \eqref{r-1} we have
\begin{eqnarray*}
&&\underset{r\downarrow0}{\lim}\frac{\log \Pi_*\eta (B(x,r))}{\log r}\\
&\geq&
\underset{q\to\infty}{\lim}\frac{\frac{q(q+1)}{2}h(\mu,\sigma)-\log 3+q \log(1-\delta)-
\underset{i=1}{\overset{q}{\sum}}i\epsilon_i}{\frac{(q+1)(q+2)}{2}\lambda(\mu,\sigma) +
5\sum_{i=1}^{q+1}ik_i\epsilon_i}\\
&=& \frac{h(\mu,\sigma)}{\lambda(\mu,\sigma)}.
\end{eqnarray*}
 Then the  desired result follows. \hfill $\Box$

\section{Proof for Corollary \ref{main-3}}\label{coro-1}

We only need to show that $D(\alpha)$ is continous in ${\rm ri}(\LL_F)\setminus A.$
At first we show that $D(\alpha)$ is upper semi-continuous on $\mathcal{L}_{F}\setminus A$. Fix $\alpha\in \mathcal{L}_{F}\setminus A$.
 Assume $\alpha_{n}\in \mathcal{L}_{F}\setminus A$, and $\lim_{n\rightarrow\infty}\alpha_n=\alpha$. By theorem \ref{main-2}, for any $\epsilon>0$, there exists $\mu_n\in \mathcal{M}(\Lambda,T)$ such that
\[\int F d\mu_{n}=\alpha_n \quad \text{ and }\quad  D(\alpha_n)<\frac{h(\mu_n,T)}{\lambda(\mu_n,T)}+\epsilon.\]
Assume $\mu_{n}$ converges to $\mu$ weakly. Recall that  the metrical entropy is upper semi-continuous, thus $\limsup_{n\rightarrow\infty}h(\mu_n,\sigma)\leq h(\mu,\sigma).$ We also have
$$
\lim_{n\rightarrow\infty}\lambda(\mu_n,\sigma)=\lambda(\mu,\sigma)\ \ \text{ and }\ \
\int fd\mu_n=\lim_{n\rightarrow\infty}\int f d\mu=\alpha,\\
$$
 thus we have
\[\limsup_{n\rightarrow\infty} D(\alpha_n)\leq \frac{h(\mu,\sigma)}{\lambda(\mu,\sigma)}+\epsilon\leq D(\alpha)+\epsilon.\] Since $\epsilon$ is arbitrary,   $D(\alpha)$ is upper semi-continuous.

%To prove the lower semi-continuity of $D(\alpha)$ we need the following simple fact, the proof of which is basic and will be omitted.
%\begin{lemma}
%Assume $V$ is a vector space, $\alpha_1,\dots \alpha_m$ are $m$ linear independent vector in $V$. If there exists sequence $\{k_{j}^{i}\},  1\leq i\leq m$, and $1\leq j<\infty$. such that
%$\lim\limits_{j\rightarrow \infty}\sum_{i=1}^{m}k_{j}^{i}\alpha_{i}=0,$
%then $\lim\limits_{j\rightarrow\infty}k_{j}^{i}=0$, for all $i$.
%\end{lemma}

Now we assume  $\alpha\in \ri (\mathcal{L}_{F}) \setminus A$. Since $A$ is closed, $\ri (\mathcal{L}_{F}) \setminus A$ is a relative open set in $\LL_F.$ Assume $\LL_F$ has dimension $l.$ Thus  there exists an $l$ dimensional simplex  $\Delta \subset ri(L_{F})\setminus A$ such that $\alpha$ is at the center of $\Delta$.
Assume $\Delta={\rm Co}\{x_0,\cdots, x_l\}$.  Define
$\Delta_i= {\rm Co}\{\alpha,x_0,\cdots, x_{i-1}, x_{i+1},\cdots,x_l\}.$
Then $\{\Delta_0,\cdots, \Delta_l\}$ form a partition of $\Delta.$

Take $\alpha_n\to\alpha$, without loss of generality, we assume $\alpha_n\in \Delta$ for all $n\ge 1$.
Assume $\alpha_n\in \Delta_{i_n}={\rm Co}\{\alpha, x_{n,1},\cdots, x_{n,l}\}$, then
$$
\alpha_n=t_{n,0}\alpha+\sum_{j=1}^l t_{n,j} x_{n,j}.
$$
 Since $\lim\limits_{n\rightarrow \infty}\alpha_n=\alpha$ we have  $\lim\limits_{n\rightarrow\infty}t_{n,0}=1$ and $\lim\limits_{n\rightarrow\infty}t_{n,j}=0$ for $j=1,\cdots,n$.

 Now fix any $\mu$ such that $\int F d\mu=\alpha$ and $\lambda(\mu,T)>0.$ For each $x_{n,j}$ find
$\mu_{n,j}$ such that
$\int f d\mu_{n,j}=x_{n,j}$ and $\lambda(\mu_{n,j},T)>0.$
Define $\mu_{n}=t_{n,0}\mu+\sum_{j=1}^{l}t_{n,j}\mu_{n,j}$. thus we have
\begin{align*}
\int F d\mu_n=& t_{n,0}\alpha+\sum_{j=1}^l t_{n,j} x_{n,j}=\alpha_n\\
 h(\mu_{n},T)=&t_{n,0}h(\mu,T)+\sum_{j=1}^{l}t_{n,j}h(\mu_{n,j},T)\\
 \lambda(\mu_{n},T)=&t_{n,0}\lambda(\mu,T)+\sum_{j=1}^{l}t_{n,j}\lambda(\mu_{n,j},T).
\end{align*}
Since $\lim\limits_{n\rightarrow\infty}t_{n,0}=1$ and $\lim\limits_{n\rightarrow\infty}t_{n,j}=0$ for $j=1,\cdots,n$,  we have
$$
\lim\limits_{n\rightarrow\infty}h(\mu_n,T)=h(\mu,T),\ \ \ \lim\limits_{n\rightarrow\infty}\lambda(\mu_n,T)\rightarrow \lambda(\mu,T).
$$
Thus by Theorem \ref{main-2} we have
 $$
 \frac{h(\mu,T)}{\lambda(\mu,T)}=\lim\limits_{n\rightarrow\infty}\frac{h(\mu_n,T)}{\lambda(\mu_n,T)}\leq \liminf_{n\rightarrow\infty} D(\alpha_n).
 $$
Now use Theorem \ref{main-2} again we conclude that $D(\alpha)\le\liminf_{n\rightarrow\infty} D(\alpha_n). $Thus  $D(\alpha)$ is lower semi-continuous.

%\section{An example}\label{example}
%
%Here we state a famous example of non-uniformly hyperbolic map.
%\begin{example}
%Let $T:[0,1]\rightarrow [0,1]$ be the Manneville-Pomeau map defined by $Tx=x+x^{1+\beta} \mod 1$, where $0<\beta<1$.  Let $\phi:[0,1]\rightarrow R^{m}$ be continuous map, take it as $\phi=(\phi_1,\phi_2,\dots,\phi_m)$.
%
%For any $\alpha \in \mathcal{L}_{\phi}\backslash{\phi(0)}$, it follows from {\bf Theorem 2} we have
%\[\dim_{H}\Lambda_{\alpha}=\sup_{\mu\in \mathcal{M}(T,\Lambda)}\left\{\frac{h(\mu,T)}{\lambda(\mu,T)},\int f d\mu=\alpha, \lambda(\mu,T)>0\right\},\]
%and  $\dim_{H}\Lambda_{\phi(0)}=1$.
%\end{example}
%

\end{document}